\documentclass[11pt,a4paper]{article}
\usepackage{epsfig,amsmath,amssymb,enumerate,mathrsfs,theorem}

\newtheorem{thm}{Theorem}[section]

\newtheorem{lem}[thm]{Lemma}

\def\pn{\par\noindent}

\title{The Rational Distance Problem for Isosceles Triangles with one rational side}
\author{Roy Barbara \and Antoine Karam }
\date{}
\begin{document}

\maketitle
\begin{abstract}
For a triangle $\Delta$ , let (P) denote the problem of the
existence of points in the plane of $\Delta$,  that are at
rational distance to the 3 vertices of $\Delta$. Answer to (P) is
known to be positive in the following situation:\\ $\Delta$ has
one rational side and the square of all sides are rational.\\
Further, the set of solution-points is dense in the plane of
$\Delta$.   See [3] \\ The reader can convince himself that the
rationality of one side is a reasonable minimum condition to set
out, otherwise problem (P) would stay somewhat hazy and scattered.
Now, even with the assumption of one rational side, problem (P)
stays hard. In this note, we restrict our attention to isosceles
triangles, and provide a \textit{complete description} of such
triangles for which (P) has a positive answer.
\end{abstract}

\baselineskip13.5pt
\section{The results}

An isosceles triangle with one rational side has one of the forms $(\lambda,\theta,\lambda)$ or $( \theta,2\lambda,\theta)$ with $\lambda \in \mathbb{Q}$ and $\theta\,\in \mathbb{R}$. Since problem (P) is invariant by a rational re-scaling, it suffices to focus on triangles of one of the forms $(1,\theta,1)\;\;\; $or$\;\;\;(\theta,2,\theta)$.\\* Note that for triangles $(1,\theta,1)$, the apex will always be disregarded and problem (P) in this case rather asks for points other than the apex.\\
Our main results are:\\

\begin{thm}Let $\Delta=(1,\theta,1)$ be an isosceles triangle with \;$\theta \in \mathbb{R}$,
$ 0\leq \theta \leq 2$. Then, answer to (P) is \emph{positive} if
and only if $ \theta^2$ has the form
 \begin{center}
 $\theta^2=2(1+pq\pm\sqrt{(1-p^2)(1-q^2)})$
 \end{center}
for some rational numbers p,q with $-1\leq p,q \leq+1$\\
\end{thm}
\begin{thm}Let $\Delta=(\theta,2,\theta)$ be an isosceles triangle with $\theta \in \mathbb{R},\; \; \theta \geq 1.$ Then, answer to (P) is \emph{positive} if and only if the main altitude $\Phi=\sqrt{\theta^2-1}$ has the form\\
\begin{center}
  $\Phi=\pm\sqrt{(p^2-1)(1-q^2)}\pm\sqrt{r^2-p^2q^2}$\\
\end{center}
for some rational numbers p, q, r with $p\geq 1 \geq q \geq 0,\;\; r \geq pq $, and where the right member is nonnegative.
\end{thm}
Regarding either triangles  $(1,\theta,1)$ \; or triangles  $(\theta,2,\theta)$,\; we say that $\theta$ \;is "suitable" if\\answer to (P) is positive for the triangle $(1,\theta,1)$, respectively $(\theta,2,\theta)$.\\
Here are the first properties or consequences of theorems 1.1 and 1.2, that will be enlightened in section 4.\\

\begin{itemize}
  \item  The suitable real numbers $\theta$ are algebraic numbers of degree at most 4.
  \item  Regarding triangles $(1,\theta,1)$, an effective procedure allows to deciding whether a given algebraic real number (of degree $\leq 4$) is suitable or not.
  \item  Given a suitable $\theta$, an effective procedure allows to construct one (possibly more) solution-point to (P).
  \item  In \textit{contrast} with the result in [3], when (P) has a positive answer, the set of solution-points in not in general dense in the plane of
  $\Delta$. More precisely, when $\theta^2$  is
  \textit{irrational}, the solution-points lie all on the union of
  two lines.
\end{itemize}

\section{proof of theorem 1}

Properties Q1 and Q2 are easily checked.
\begin{itemize}
  \item (Q1) Set $Z=\{z \in \mathbb{R},\;z=pq\pm\sqrt{(1-p^2)(1-q^2)},\; p, q \in \mathbb{Q},\; -1\leq p,q \leq+1\}$. Then, \begin{enumerate}
                                                                                                                    \item For $z \in Z$\;\; we have\;\; $-1\leq z \leq +1$.
                                                                                                                    \item $Z$ is closed by opposite $[z\in \mathbb{Z} \Rightarrow -z \in \mathbb{Z}].$
                                                                                                                  \end{enumerate}
  \item (Q2) Let $\Delta = ABC$ be a triangle with $AB=AC=1.$ Let M be a point in the plane of $\Delta,\; M\neq A.$  Set $R=MA,\, S=MB,\, T=MC$. Define $u=\frac{1}{2}(R^2-S^2+1)$ and $v= \frac{1}{2}(R^2-T^2+1)$.\; Then,\; $ u^2 \leq R^2 \; \; \;$ and$\; \; \; v^2\leq R^2.$
\end{itemize}
We also need 3 lemmas:
\begin{lem}
Let $\Delta = ABC$ be a triangle with $AB=AC=1$ and $BC=\theta$ $(0\leq \theta \leq 2)$. Suppose that $\theta^2 \in \mathbb{Q}$. Then, there are (infinitely many) points M in the plane of $\Delta, \, M\neq A$, such that MA, MB, and MC are all rational numbers.
\end{lem}

\textbf{Proof:} Set $w=\angle BAC , a =\cos w, b = \sin w.$ By the law of cosines, $a=1- \frac{\theta^2}{2}$. Since $\theta^2 \in \mathbb{Q}$, then a $\in\mathbb{Q}$. For $ \psi \in \mathbb{Q}-\{a,\pm 1\}$, set $x=\frac{\psi^2-1}{2(\psi-a)}$. Then, $x \in \mathbb{Q}-\{0\}$. Let M be on the axis $\overrightarrow{AB}$ with $\overline{AM}=x$. As $ x \neq 0, $ then $M \neq A$. Since $MA = |x| $ and $MB=|x-1|$, then, $MA, MB \in \mathbb{Q}.$\;\; Now, using Pythagoras, we may write:\\
$\overline{MC}^2=(x-a)^2+b^2=x^2 -2ax+a^2+b^2 = x^2-2ax+1 \\ =$
$\frac{\psi^4-4a\psi^3+4a^2\psi^2+2\psi^2-4a\psi+1}{4(\psi-a)^2} = (\frac{\psi^2-2a\psi+1}{2(\psi-a)})^2$. Hence, $MC \in \mathbb{Q}$.  \maltese\\
\\
\begin{lem}
Let $p, q \in \mathbb{Q}$ with $p+q\neq0$ and $|p-q| \neq 2$.
Then, there are rational numbers $R,S,T, \; \; R \neq 0,$   such that\\
\begin{center}$\frac{R^2-S^2+1}{2R}=p\;\;\;$\;\;\;\;\;         and \;\;\;\;\;       $\;\;\;\frac{R^2-T^2+1}{2R}=q$\\ \end{center}
\end{lem}
 \textbf{Proof:}\\The following values will do:\;\;\; $R=\frac{4-(p-q)^2}{4(p+q)}.$\;\;\;\; $\; S=\frac{4+2pq-3p^2+q^2}{4(p+q)}$.\;\;\;\; $T= \frac{4+2pq+p^2-3q^2}{4(p+q)}.\;\;\maltese$
\newline
\begin{lem}
Let $\Delta =ABC$ be a non-degenerated isosceles triangle with $AB=AC=1$ and $w=\angle BAC, \; \; 0 < w <\pi$.\;\;  Set $a = \cos w$.\\
Then, the following statements are equivalent:
\begin{description}
  \item[(i)] There is a point M in the plane of $\Delta , \; M\neq A,$\; such that MA, MB, MC are all rational numbers.
  \item[(ii)] There are rational numbers R, S, T, $\; R \neq 0$,\; such that,\; if \; $u= \frac{R^2-S^2+1}{2}$ \;\; and  \\ $v= \frac{R^2-T^2+1}{2}$,\;\; we have
    \begin{equation}
      R^2a^2+u^2+v^2 = R^2+2auv
    \end{equation}
\end{description}
\end{lem}
\textbf{proof:}
Set $b=\sin w > 0$. Consider a x-y system such that A(0,0), B(1,0), and C(a,b).\\
$(i)\Rightarrow (ii):$ Let $M(u,\rho)$ be a solution-point to (P), $M\neq A$. Set $R=MA, \;\; S=MB,\\ T=MC$. Then, $R,S,T \in \mathbb{Q},\;\; R>0$.\;\; By Pythagoras we may write:\\
\begin{center}
$u^2+\rho^2=R^2$\\
$(u-1)^2+\rho^2=S^2$\\
$(u-a)^2+(\rho-b)^2=T^2$
\end{center}
From the first two relations we get $2u=R^2-S^2+1$. Set also
$2v=R^2-T^2+1$. With $u^2+\rho^2=R^2$, $\;a^2+b^2=1,$ the third
equation yields $R^2+1-2ua-2\rho b=T^2$. Hence, $2ua+2\rho b=2v$,
\;so\;
$\rho b = v-ua. $\;\; Therefore, $\rho^2 b^2=(v-ua)^2$, that is, $(R^2-u^2)(1-a^2)\\ =(v-ua)^2$.\; Rearranging, we get (1). \\
$(ii)\Rightarrow (i)$: Suppose (1) satisfied with some $R, S, T \in \mathbb{Q}, \; R\neq 0$, and with $u=\frac{(R^2-S^2+1)}{2}$,  $v=\frac{(R^2-T^2+1)}{2}$.\; Rewrite (1) as $(R^2-u^2)(1-a^2)=(v-ua)^2$,\;  that is, \\
\begin{equation}
(R^2-u^2)b^2=(v-ua)^2.
\end{equation}
Since $b^2 >0$ and the right member in (2) is nonnegative, we get $R^2-u^2 \geq 0$. Define $\rho_0 = \sqrt{R^2-u^2} \in \mathbb{R}^+$. (2) becomes $(\rho_0b)^2= (v-ua)^2$, hence\; $\pm \rho_0b =v-ua$.\;  Let $\rho \in \{\pm \rho_0\}$ such that  $\rho b=v-ua $.\; We then have:\\
\begin{eqnarray}
& &u^2+\rho^2 = R^2\\
& &\rho b = v-ua
\end{eqnarray}

Consider the point $M(u,\rho)$. Since $MA^2=R^2 >0$, then, $M\neq A$ and $MA \in \mathbb{Q}$. Since $MB^2=(u-1)^2+\rho^2=R^2+1-2u=S^2,$ then, $MB\, \in \mathbb{Q}.$ Finally, from $MC^2 = (u-a)^2+(\rho-b)^2 = R^2+1-2ua-2\rho b=R^2+1-2ua-2(v-ua)=R^2+1-2v=T^2$, we get $ MC \in \mathbb{Q}$.  \maltese \\
Note that a related characterization is to appear in [2].\\
\\
* We now are ready to prove theorem 1:
\\
Let $\Delta = ABC$ be a triangle with $AB=AC=1$ and $BC = \theta $, $\; 0\le \theta \le 2$.\\
Set $w=\angle BAC,\;\; a=\cos w,\;\; b=\sin w$.\\
If $\Delta$ is degenerated $(\theta=0$ or $2)$, as quickly seen, both parts (i) and (ii) hold.\\
From now on, we assume $\Delta$ non-degenerated. Thus, $0 < \theta < 2, \; 0<w< \pi, \; $ and $b=\sin w >0$.\\
\\
$(i)\Rightarrow(ii)$: Assume that (P) has a positive answer. By lemma 2.3, there are R, S, T $\in \mathbb{Q},\; R\neq 0$,\;\; such that relation (1) holds with\; $u=\frac{R^2-S^2+1}{2}$\;\; and\;\; $v=\frac{R^2-T^2+1}{2}$. \\
It follows that a is a zero of the trinomial in t\\ \begin{center} $R^2 t^2-2uvt+(u^2+v^2-R^2)=0$.\\ \end{center}
Hence, $ a=(\frac{u}{R})(\frac{v}{R})\pm \sqrt{(1-(\frac{u}{R})^2)(1-(\frac{v}{R})^2)}$. Set $p=(\frac{u}{R})$ and $q=(\frac{v}{R})$. Then $p,q \in \mathbb{Q}$. According to (Q2), $u^2\leq R^2,\; v^2 \leq R^2$,  that is, $p^2,\; q^2\leq1$, where $a=\cos w= pq\pm \sqrt{(1-p^2)(1-q^2)}$. By (Q1), -a has the same form than a. For convenience, we rather put\\
\begin{center}
$-a=-\cos w = pq \pm \sqrt{(1-p^2)(1-q^2)}$
\end{center}
By the law of cosines,\; $\theta^2=2(1-a)$.\;\; Hence,\\
\begin{center}
$\theta^2 = 2(1+pq \pm \sqrt{(1-p^2)(1-q^2)})\;$\;\; where $p,q
\in \mathbb{Q},\; \; -1\leq p,q \leq +1$
\end{center}
$(ii)\Rightarrow(i)$: Suppose that $\theta^2 = 2(1+pq \pm \sqrt{(1-p^2)(1-q^2)})$ for some $p, q \in \mathbb{Q}$, $-1\leq p,q \leq +1$. Since $\theta^2 =2(1-a)$, we get $-a= pq \pm \sqrt{(1-p^2)(1-q^2)}$. By (Q1), $-(-a)$ has the same form than -a.\;   For convenience, we rather set:\\
\begin{equation}
a=\cos w= pq \pm \sqrt{(1-p^2)(1-q^2)}
\end{equation}
where $p,q \in \mathbb{Q},\; -1\leq p,q \leq +1.$\;\; Therefore $(a-pq)^2=(1-p^2)(1-q^2)$ and hence \\
\begin{equation}
a^2+p^2+q^2=1+2apq
\end{equation}
\begin{itemize}
  \item If $p+q=0,$ then $p^2=q^2$,\; so \; $a=-p^2\pm \sqrt{(1-p^2)^2}\;=-p^2\pm(1-p^2)$.\\ Hence $a=\cos w \in \mathbb{Q}$.\; Hence $\theta^2=2(1-a) \in \mathbb{Q}$.\; Lemma 2.1 gives the result.
  \item We assume now $p+q\neq 0$. Claim $|p-q| \neq 2$: otherwise, if $|p-q|=2$ and since $p,q \in [-1,+1]$, we would get $\{p,q\}=\{\pm1\}$,\; and hence by (5), $ a= \cos w =-1$,\\ that is \;$w=\pi$,\; a contradiction.
\end{itemize}
Now, we apply lemma 2.2 . There are $R, S, T \in \mathbb{Q}, \; R \neq 0$, such that:\\
\begin{center}
  $p=\frac{R^2-S^2+1}{2R}\;\;\;$ \;\;\;\;\;   and \;\;\;\;\;    $\;\;\; q=\frac{R^2-T^2+1}{2R}$
\end{center}
Set $ u=\frac{1}{2}(R^2-S^2+1)\;$ and $ v=\frac{1}{2}(R^2-T^2+1).
\; $ We have \; $p=\frac{u}{R}$ \; and \; $q=\frac{v}{R}$.
Replacing p and q respectively by $\frac{u}{R}$ and $\frac{v}{R}$
in (6) yields
\begin{center}
  $R^2a^2+u^2+v^2=R^2+2auv$
\end{center}
which is relation (1). \; Lemma 2.3 achieves the proof.   \maltese
\\
\section{Proof of theorem 2}
We need two lemmas:
\begin{lem}
Let $\Delta =ABC$ be a triangle with $AB=AC=\theta \;$ and $\; BC=2$ \\
$(\theta \in \mathbb{R}, \; \theta \geq 1)$. Suppose that $\theta^2 \in \mathbb{Q}$. Then,\\
\begin{description}
  \item[(i)] There are (infinitely many) points in the plane of $\Delta$, that are at rational distance to the 3 vertices of $\Delta$.
  \item[(ii)] The main altitude $\Phi=\sqrt{\theta^2-1}$ can be put in the form\\
  $\Phi=\sqrt{(p^2-1)(1-q^2)}+\sqrt{r^2-p^2q^2}$, with $p, q, r \in \mathbb{Q}$, $p\geq 1\geq q \geq 0, \; r \geq pq$.\\
\end{description}
\end{lem}
\textbf{Proof}: Let O be the midpoint of BC.
\begin{description}
  \item[case1]: T is degenerated $(\theta=1, \Phi=0)$: Both parts (i) and (ii) are obvious.
  \item[case2]: T is non-degenerated $(\theta >1, \Phi >0):$ Since $\Phi^2=\theta^2-1 \in \mathbb{Q}$ and $\Phi >0$, set $\phi^2=f, \; f\in \mathbb{Q}, f>0.$ Select a positive integer N such that $Nf\geq 2$ (infinite choice). For (i), let M be one of the two points on BC such that $MO = Nf-\frac{1}{4N}$. Then, $MB,MC \in \mathbb{Q}\;$ and $MA^2=MO^2+OA^2=MO^2+\Phi^2=(Nf-\frac{1}{4N})^2+f=(Nf+\frac{1}{4N})^2$, so $MA \in \mathbb{Q}$. For (ii), the values $p=Nf-\frac{1}{4N}$, $q=1$, and $r=Nf+\frac{1}{4N}$ will do.\;\;  \maltese
\end{description}
\begin{lem}
Let $e \in \mathbb{Q}$ and $\alpha ,\beta \in \mathbb{Q}$ -\{0\}.\; Suppose that\\
\begin{equation}
(\frac{e}{\beta}-\beta) = (\frac{e}{\alpha}-\alpha)\, +\, 4
\end{equation}
Then, for some $p, \, q \in \mathbb{Q},\; p\neq \pm1,\; q \neq 1 \; (if \, e \neq 0, \; q\neq \pm 1)$ such that\\
\begin{center}
$\alpha +\beta=\pm 2p(1-q)$\\
$\;\; \alpha \beta =(p^2-1)(1-q)^2$\\
$e=(p^2-1)(1-q^2)$
\end{center}
\end{lem}
\textbf{Proof}: (7) clearly shows that $\alpha \neq \beta $. \;
Set $\gamma = \alpha \beta$ \; and \; $\delta = \alpha - \beta$.
\\Then, $\gamma,\; \delta \in \mathbb{Q}-\{0\}$. From (7) we get $e(\frac{\alpha - \beta}{\alpha \beta})=4-(\alpha -\beta)$. That is, $e\frac{\delta}{\gamma}=4-\delta$. Hence,\\
\begin{equation}
e=\gamma (\frac{4-\delta}{\delta})
\end{equation}
Now, $\alpha$ and $-\beta$ are the roots of \;\; $t^2-\delta
t-\gamma =0$. \; Therefore, the discriminant must be \\
a rational square,\;  say $\delta^2+4\gamma=\epsilon^2,\; \;
\epsilon \in \mathbb{Q}$.
\\ We have
\begin{equation}
\gamma=\frac{\epsilon^2-\delta^2}{4}
\end{equation}
$\alpha$ and $-\beta$ are in some order $\frac{\delta+\epsilon}{2}$ and $\frac{\delta-\epsilon}{2}$. In all cases, we have
\\
\begin{equation}
\alpha +\beta=\pm \epsilon
\end{equation}
Set\; $p=\frac{\epsilon}{\delta}$\;  and \; $ q=
1-\frac{\delta}{2}$.\; Then:
\begin{equation}
\delta=2(1-q)
\end{equation}
and
\begin{equation}
\epsilon = 2p(1-q)
\end{equation}
From (11) we easily get
\\
\begin{equation}
\frac{4-\delta}{\delta} = \frac{1+q}{1-q}
\end{equation}
\\
Finally,
\\ \begin{itemize}
    \item From (10) and (12) we obtain \; $\alpha + \beta =\pm \,2p(1-q).$
    \item Using (9), (12), and (11), we may write
    \\$\alpha \,\beta=\gamma=\frac{\epsilon^2-\delta^2}{4} = \frac{4p^2(1-q)^2-4(1-q)^2}{4} = (p^2-1)(1-q)^2$.
    \item From (8), this latter, and (13), we may write
    \\$ e = \gamma (\frac{4-\delta}{\delta}) = (p^2-1)(1-q)^2 \, \frac{1+q}{1-q} = (p^2-1)(1-q^2)$.
    \maltese\\
  \end{itemize}
* We now are ready to prove theorem 2:\\
Let $\Delta$ = ABC be a triangle with $AB=AC=\theta$\;  and
$BC=2\; (\theta \in \mathbb{R},\; \theta \geq 1)$. Let $\Phi =
\sqrt{\theta^2 -1}$ be the main altitude. Let O be the midpoint of
BC. Consider the x-y axes with origin O, where
$\overrightarrow{OC}$ defines the x-axis and $\overrightarrow{OA}
$ the y-axis. We have the coordinates:
\\$A(0,\Phi),\; B(-1,0)$,\; and \; $C(1,0)$.\\

(ii)\;$\Longrightarrow$\;(i) \;Suppose that\; $\Phi = \mu + \nu, \;\; \mu = \pm \, \sqrt{(p^2-1)(1-q^2)},\;\; \nu = \pm \, \sqrt{r^2-p^2q^2}, $\;\;  as in (ii).
Consider any of the points $M(\pm pq,\,\;\mu)$. We may write\\
$MA^2 =p^2q^2+(\Phi-\mu)^2 = p^2q^2 + \nu^2 = r^2, \; \; MB^2 = (\pm pq+1)^2 +\mu^2$ $=$\\
$p^2q^2 \pm 2pq+1+(p^2-1)(1-q^2)=p^2\pm 2 pq + q^2 = (p \pm q)^2$,\; and similarly $ MC^2 = (p \mp q)^2$.\\ It follows that MA, MB and MC are \textit{all} rational distances.\\

(i) $\Longrightarrow$ (ii) Suppose that some point $M(x_0,y_0)$
lying in the plane of $\Delta$ satisfies $MB = R \in \mathbb{Q}$,
$ MC = S \in \mathbb{Q}$, and $MA = r \in \mathbb{Q}$. Set $
\underline{y_0^2 = e}$ and $\underline{(\Phi - y_0)^2 = f }$. The
pythagorean relations are
\\
\begin{eqnarray}
(x_0 + 1 )^2 + e = R^2\\
%
%
(x_0 - 1 )^2 + e = S^2\\
%
%
x_0^2 + f = r^2
\end{eqnarray}
Subtracting (14) and (15) gives $4x_0 = R^2 - S^2$.\;  Hence, $\underline{x_0 \in \mathbb{Q}}$.  From $x_0 \in \mathbb{Q}$, (14) and (16), we get $ e, \, f \in \mathbb{Q}$.\; Hence, $y_0 = \pm \sqrt{e} , \,\; \Phi - y_0 = \pm \sqrt{f}$,\; with $e, f \in \mathbb{Q}^+$. In particular,\\
\begin{equation}
\Phi\;=\; y_0\pm\sqrt{f}\;=\;\pm \sqrt{e} \pm \sqrt{f}
\end{equation}
\begin{itemize}
  \item If $e = 0$, then $\Phi^2 = f$, so $\theta^2 = \Phi^2-1 = f-1 \in \mathbb{Q}$. In this case, lemma 3.1 gives the result.
  \item From now on, we assume that $\underline{e > 0}$.\\
  \end{itemize}
  Rewrite (14) and (15) as :
  \\
  \begin{eqnarray}
  e = (R-(x_0+1))(R+(x_0+1))\\
%
%
%
  e = (S-(x_0-1))(S+(x_0-1))
  \end{eqnarray}
  \\Set $\alpha =S-(x_0-1)$ and $\beta = R-(x_0+1)$. Clearly, as $e \neq 0, \, \alpha , \, \beta \in \mathbb{Q}-\{0\}$.
  \\Subtracting $R+(x_0+1) = \frac{e}{\beta}$\; and $R-(x_0+1)=\beta$ \; yields  \; $2x_0+2 = \frac{e}{\beta}-\beta$.
  \\Similarly, $S+(x_0-1) = \frac{e}{\alpha}$\;\; and \; $S-(x_0-1)=\alpha\;$  yield \;  $2x_0-2 = \frac{e}{\alpha}-\alpha$.
  \\Summing and subtracting these two relations provides
  \\
  \begin{equation}
  4x_0 = (\alpha + \beta )(\frac{e-\alpha \beta}{\alpha \beta})
  \end{equation}
\begin{equation}
  (\frac{e}{\beta}-\beta)=(\frac{e}{\alpha}-\alpha)+4
\end{equation}
From (21) and lemma 3.2, we deduce the existence of $p,\,q\in \mathbb{Q}, \; p,q \neq \pm 1$, such that
\\
\begin{center}
$\alpha +\beta=\pm 2p(1-q)$,\;\;\;\;\; $\alpha \beta
=(p^2-1)(1-q)^2$,\;\;\;\;\; $e=(p^2-1)(1-q^2)$.
\end{center}
Using this and (20), we may write:
\\$4x_0=\pm 2p(1-q)(\frac{(p^2-1)(1-q^2)-(p^2-1)(1-q)^2}{(p^2-1)(1-q)^2})= \pm \, 4pq$. \; Hence, $x_0 = \pm \, pq$ \\
Now, by (16),$\;\;f = r^2 - x_0^2 = r^2 - p^2q^2$,\;\;and
consequently
\\
\begin{center}
   $\Phi = \pm \, \sqrt{e} \pm \sqrt{f} = \pm \, \sqrt{(p^2-1)(1-q^2)} \pm \, \sqrt{r^2-p^2q^2}$.
\end{center}
Finally and without loss of generality, in such expression of
$\Phi$, we may assume that p, q, r are $\textit{nonnegative}$.
From $f \geq 0$, we get $r^2\geq p^2q^2$, so $ r\geq pq$. From $ e
= (p^2-1)(1-q^2) > 0$, we see that $(p^2-1)$ and $(q^2-1)$ have
opposite signs. Up to a permutation of p and q, we always may
assume that $p^2-1 > 0 \;$,\; hence $p^2 > 1 > q^2 \;$,\; so
$p>1>q$.  \maltese
\section{First Consequences}
\begin{description}
\item[$\blacksquare$] Note that in both theorems 1.1 and 1.2, a suitable $\theta$ satisfies $\theta^2=\mu\pm\sqrt{\nu}$, $\mu,\,\nu \in \mathbb{Q},$\\$\nu \geq 0.$ \;Hence $\theta$ must be an algebraic number of degree 1, 2, or 4. In particular:
    \begin{description}
    \item[-] If $\theta$  is transcendental or has algebraic degree 3 or $\geq 5,\; \; \theta$  is \textit{not} suitable.
    \item[-] If $\theta$  has algebraic degree 4, whence $\theta^2$  has also degree 4 \\(ex. $\theta=\frac{1}{4}(1+\sqrt{3}+\sqrt{5}))$, then $\theta$  is
\textit{not} suitable.
    \item[-] If $\theta^2 \in \mathbb{Q}$,  then $\theta$  is always suitable.\\
    \end{description}
    We focus now on the class of algebraic numbers $\theta$  of degree 2 or 4 satisfying
    \begin{center}
    $\theta^2=a \pm\sqrt{b}$,\;  $a, b \in \mathbb{Q}, \; \; b>0,\;\; \sqrt{b} \notin \mathbb{Q}$
    \end{center}
    Theorems 1.1 and 1.2 give a satisfactory answer. Moreover, one may ask whether a given real number in this class can be recognized as suitable or not by an effective procedure? At least regarding triangles $(1,\,\theta,\,1)$, we answer now positively:\\
\item[$\blacksquare$] Let $\theta \, \in \mathbb{R},\; 0 < \theta <2$, be given, where $\theta^2=a\pm \sqrt{b},\; a,\;b \in \mathbb{Q},\; b>0,\; \sqrt{b} \notin \mathbb{Q}.$\\
    In field theory one shows that if $\theta^2=c\pm \sqrt{d},\; c,\;d \in \mathbb{Q},\; d\geq 0$,  then,  $c=a$\\and  $d=b$. Therefore, assuming that $\theta^2=(2+2pq)\pm \sqrt{4(1-p^2)(1-q^2)}$  as in theorem 1.1 would imply  $2+2pq=a$ \; and \; $4(1-p^2)(1-q^2)=b$,\;  that would lead to \;$p^2q^2=\frac{(a-2)^2}{4}$ \; and \; $p^2+q^2 = \frac{(a-2)^2+4-b}{4}$. \; The algorithm is then: Find the roots $t_1$ and $t_2$\; of
    \begin{center}
    $f(t)=t^2-(\frac{(a-2)^2+4-b}{4})t+\frac{(a-2)^2}{4}=0$
    \end{center}
    If $t_1,\; t_2$ lie in $\mathbb{Q}\bigcap[0,1]$ and if $t_1$ and $t_2$ are both rational \textit{squares}, then $\theta$ is suitable, otherwise $\theta$ is not.
\item[$\blacksquare$] There is an effective procedure to finding
solution-points when $\theta$ is suitable. Regarding triangles
$(1,\;\theta,\;1)$, such algorithm can be extrapolated from lemmas
2.1, 2.2, 2.3, theorem 1.1, and their proofs.  Regarding triangles
$(\theta,\;2,\;\theta)$,  this is immediate: If $\Phi =
\sqrt{\theta^2-1} =
\epsilon\sqrt{(p^2-1)(1-q^2)}+\epsilon'\sqrt{r^2-p^2q^2}$,
$\epsilon,\; \epsilon' \in\{\pm1\}$,  as in theorem 1.2,
solution-points are
    \begin{center}
    $M\big(\pm pq,\;\; \;\epsilon\sqrt{(p^2-1)(1-q^2)}\;\big)$\\
    \end{center}
\item[$\blacksquare$] Finally,we show that the set of solution-points is not in general dense in the plane of the triangle. More precisely we prove the following when $\theta^2$\; is \textit{irrational}:\\
- If $\Delta=(\theta,\;2,\;\theta)$,\;all solution-points lie on
the union of 2 lines that are parallel to the basis of
 $\Delta$.\\
- If $\Delta =(1,\;\theta,\;1)$,\;all solution-points lie on the
union of 2 concurrent lines at the apex,\; that are symmetric
through the main altitude.\\
$ $\\
$\bigstar$ Let $\Delta =( \theta,\;2,\;\theta),\; \theta > 1$,
where $\Phi=\sqrt{\theta^2-1}
=\epsilon\sqrt{a}+\epsilon'\sqrt{b},\;\epsilon, \epsilon'
\in\{\pm1\},\; a,\; b \in\mathbb{Q}, \;a,\, b >0$. We assume that,
either $\sqrt{a}$  and $\sqrt{b}$  are non-degenerated and
non-associated radicals $(\sqrt{ab} \notin \mathbb{Q})$, or, that
exactly one of $\sqrt{a},\; \sqrt{b}$ is degenerated. This (most
frequent) situation corresponds precisly to the fact that
$\theta^2$  is irrational. In field theory, one then proves the
following: If $\Phi = \eta \sqrt{c}+\eta'\sqrt{d},\;\; \eta'\,
\eta' \in\{\pm1\},\; c,\, d \in \mathbb{Q},\; c,d \geq0$,  then we
must have $(c,d,\eta,\eta')=(a,b,\epsilon,\epsilon')$ or
$(b,a,\epsilon',\epsilon)$. In particular, $\{\eta\sqrt{c},\,
\eta' \sqrt{d}\}=\{\epsilon\sqrt{a},\epsilon' \sqrt{b}\}$. Now
suppose that $\theta$ is suitable, that is,
$\Phi=\sqrt{\theta^2-1}=\eta\sqrt{(p^2-1)(1-q^2)}+\eta'
\sqrt{r^2-p^2q^2},\; \eta,\, \eta' \, \in \{\pm1\}$, as in theorem
1.2.\; By the above property, we must have
\begin{center}
$\eta\sqrt{(p^2-1)(1-q^2)}\in \{\epsilon\sqrt{a},\;\;\epsilon'
\sqrt{b}\}$
\end{center}
By the proof of theorem 1.2, any solution-point $M(x_0,y_0)$ satisfies\\
\begin{center}
$x_0=\pm pq \;$\;\;\;\;\;\;\;   and   \;\;\;\;\;\;\;$y_0=\eta
\sqrt{(p^2-1)(1-q^2)}$
\end{center}
Hence $y_0 \in \{\epsilon \sqrt{a},\; \epsilon' \sqrt{b}\}$. Therefore, all solution-points lie on the union of the 2 lines:\\
\begin{center}
$y=\epsilon\sqrt{a}\;\;\;\;\;\;\;\;\;\;\;\;\;\;\;\;\;\;\;\;\; y=\epsilon' \sqrt{b}$. \;\;\;\;\;\;\;\;\;\;\maltese\\
\end{center}
$\bigstar$ Let $\Delta\;=\;(1,\;\theta,\;1),\;\;0<\;\theta\;<2$,\;
with apex angle $\omega,\;\;a=\cos \omega$,\; and axis of symmetry
$\Gamma$. Suppose that $\theta$ is suitable whereas $\theta^2$ is
\textit{irrational}. Denote by $\Sigma$ the set of
solution-points. According to theorem 1.1,  $a=p_0q_0\pm
\sqrt{(1-p_0^2)(1-q_0^2)}\;,\;\; p_0,q_0 \in
\mathbb{Q}\cap[-1,1],$ \;where the radical is non-degenerated as
$\theta^2=2(1-a)\;\notin \mathbb{Q}$.\; Let $M(u, \rho)\,\in
\Sigma$.\;\, Set $R=MA, \;\;S=MB,\;\;T=MC,\;\; \; R, S, T \in
\mathbb{Q},\;\; R>0.$ \; By the proofs of lemma 2.3 and theorem
1.1 \big(parts $(i)\Rightarrow(ii)\big)$,  we know that
$u=\frac{1}{2}(R^2-S^2+1)$\; and that, with
$v=\frac{1}{2}(R^2-T^2+1),\;p=\frac{u}{R},$\; and $q=\frac{v}{R}$,
we have $a=pq\pm \sqrt{(1-p^2)(1-q^2)}$.  Consequently,\; $pq\pm
\sqrt{(1-p^2)(1-q^2)} = p_0q_0\pm \sqrt{(1-p_0^2)(1-q_0^2)}$.\;
Since the latter radical is non-degenerated, one proves in field
theory that\; $pq=p_0q_0$  and \;$ (1-p^2)(1-q^2) =
(1-p_0^2)(1-q_0^2)$,\;  that yields $pq =p_0q_0$ \; and \;
$p^2+q^2=p_0^2+q_0^2$.\; It is then elementary to see that
    \begin{center}
     $p \in \{\pm p_0,\;\pm q_0\}$\\
    \end{center}
    \begin{description}
      \item[\underline{case 1:  $p_0q_0 = 0$}]$\;\;\;\; p_0=q_0=0$\; is impossible since \;$\sqrt{(1-p_0^2)(1-q_0^2)} \; \notin \mathbb{Q}.$\; Without loss of generality, assume \;$p_0\neq 0$\;\, and\;\,  $q_0 = 0.$\; Then, $p \in \{0,\; \pm p_0 \}.$\\
  If \;$ p= 0,$ \;then \;$ u=pR=0,$ \; so the point \;$M(u,\;\rho) = M(0,\;\rho)$\; lies on the y-axis,  say\;$L_0$.\\
  If $\; p=\pm p_0,$\; the ratios \;$\pm \frac{\sqrt{1-p^2}}{p}$\;  can only take 2 values \;$k_1 =\frac{\sqrt{1-p_0^2}}{p_0}$\;  and $k_2 =-\frac{\sqrt{1-p_0^2}}{p_0}$. \;From $u^2 + \rho^2 = R^2$\; and\; $u=pR$,\; we get \;$\rho^2 = R^2(1-p^2)$,\; hence $ \rho = \pm R\sqrt{1-p^2} ,$\; and hence $\frac{\rho}{u}=\frac{\pm R\sqrt{1-p^2}}{pR} = \pm\frac{\sqrt{1-p^2}}{p} \in \{k_1,\; k_2\}$.\; It follows that M lies on the union of the two lines \;$L_1:\;\; y=k_1 x$\;  and \;$L_2:\;\; y = k_2 x$.\; The reader can check that one (\textit{and only one}) line, say \;$L \in \{L_1,\; L_2\}$\; is the reflexion of $L_0$ through \;$\Gamma$.\; Since \;$\Sigma$\; is closed by symmetry through \;$\Gamma$,\; we conclude that \;$\Sigma \subseteq L_0\sqcup L$.
      \item[\underline{case 2:  $p_0q_0 \neq 0$}]$\;\;\;\ $     Since \;$p \in \{\pm p_0,\; \pm q_0\}$,\;  the ratios \;$\pm \frac{\sqrt{1-p^2}}{p}$\; can only take 4 values:
\begin{center}
$k_1= \frac{\sqrt{1-p_0^2}}{p_0}$,\;\;\;\;\; $k_2= -\frac{\sqrt{1-p_0^2}}{p_0}$,\;\;\;\;\; $k_3= \frac{\sqrt{1-q_0^2}}{q_0}$,\;\;\;\;\; $k_4= -\frac{\sqrt{1-q_0^2}}{q_0}$
\end{center}
As noted above, \;$\frac{\rho}{u} = \pm \frac{\sqrt{1-p^2}}{p}$,\;
hence \;$\frac{\rho}{u} \in \{k_1,\; k_2,\;k_3,\;k_4\}$.
Therefore, M lies on \textit{at most} the union of the 4 lines \;
$L_1, \;L_2,\;L_3,\;L_4,\;$ with respective equations
\begin{center}
$y=k_1 x,\;\;\;\;\;\;\;\;\;\;y=k_2
x,\;\;\;\;\;\;\;\;\;\;y=k_3 x,\;\;\;\;\;\;\;\;\;\;y=k_4 x.$
\end{center}
Among these 4 lines (for convenience we omit the details),  only
two lines, say \;$L, L'\;(L \in \{L_1,\;L_2\},\; L' \in
\{L_3,\;L_4\})$, are symmetric through \;$\Gamma$.\; Since
$\Sigma$ is closed by symmetry through \;$\Gamma$,\; we conclude
that \;$\Sigma \subseteq L \cup L'$.\;\;\;\maltese
    \end{description}

\end{description}
\section{Related Open Problems}

Introduce the set $\Omega = \{(p^2-1)(q^2-1),\; p,\; q \in \mathbb{Q},\; p,\, q \geq0\}$.\\
It can be proved that $\; -1,\; 2,\; \frac{1}{2}\; \notin \Omega$\;  (properties related to the Fermat quartic equation\\
  $X^4-Y^4=Z^2$)\\
On the other hand, a representation of $\omega \in \Omega$  is not
\textit{in general} unique as shown in the example
\begin{center}
$-\frac{72}{25}=\big((\frac{11}{5})^2-1\big)\big((\frac{1}{2})^2-1\big)=\big(2^2-1\big)\big((\frac{1}{5})^2-1\big)$
\end{center}
Apart from $1=(0^2-1)(0^2-1)$ ,\;\;1 has infinitely many
representations since
\begin{center}
$1=\big((\frac{z}{x})^2-1\big)\big((\frac{z}{y})^2-1\big)$
\end{center}
for any pythagorean triple (x, y, z)\;\;\;\;(x, y, z are positive integers with $x^2+y^2=z^2$).\\
Questions of interest are:\\
\begin{description}
\item[\underline{P1}] Is $\Omega$ a decidable set? (i.e. is there
an effective procedure to determine whether a given rational
number lies or not in $\Omega$  ?) \item[\underline{P2}]
Disregarding 1,\; does an element in $\Omega$\; have a finite
number of representations ?\item[\underline{P3}] Which elements in
$\Omega$ do have a \textit{unique} representation (up to the order
of the factors)? \item[\underline{P4}] For which triangles $(1,\;
\theta, \;1)$, respectively $(\theta, \;2,\; \theta)$, is the set
of solution-points to problem(P)  a finite set ?   a single set ?
\item[\underline{P5}] Is there an algorithm to decide whether an
algebraic number \;$\theta \geq 1$\; (of degree $\leq 4$) is
suitable or not for the triangle $(\theta,\;2,\;\theta)$ ?
\end{description}

\bigskip
{\footnotesize \pn{\bf Roy Barbara. \quad Antoine Karam} \;
\\Lebanese University,
Faculty of Science II.\\
Fanar Campus. P.O.Box 90656. \\
Jdeidet El Metn. Lebanon.\\
{\tt Email: roybarbara.math@gmail.com \\
\qquad amkaram@ul.edu.lb}



\end{document}